\documentclass[10pt]{amsart}
\usepackage{amscd}
\usepackage{amssymb}

\newtheorem{thm}[equation]{Theorem}
\newtheorem{cor}[equation]{Corollary}
\newtheorem{prop}[equation]{Proposition}
\newtheorem{lem}[equation]{Lemma}

\theoremstyle{definition}

\newtheorem{rem}[equation]{Remark}
\newtheorem{exa}[equation]{Example}

\newtheorem{que}[equation]{Question}
\numberwithin{equation}{section}

\newcommand{\iso}{\stackrel{\simeq}{\rightarrow}}
\newcommand{\inj}{\hookrightarrow}
\newcommand{\surj}{\twoheadrightarrow}
\newcommand{\ar}{\rightarrow}
\newcommand{\opn}{\operatorname}
\newcommand{\cat}[1]{\operatorname{\mathsf{#1}}}
\newcommand{\bdot}{{\textstyle \cdot}}

\newcommand{\mfrak}[1]{\mathfrak{#1}}
\newcommand{\mcal}[1]{\mathcal{#1}}

\newcommand{\mbf}[1]{\mathbf{#1}}
\newcommand{\mrm}[1]{\mathrm{#1}}
\newcommand{\mbb}[1]{\mathbb{#1}}

\newcommand{\tup}[1]{\textup{#1}}

\newcommand{\bwedge}{\bigwedge\nolimits}

\title[Rigid Dualizing Complex]
{The Rigid Dualizing Complex of a Universal Enveloping Algebra}
\author{Amnon Yekutieli}
\address{Department of Theoretical Mathematics,
The Weizmann Institute of Science,
Rehovot 76100, ISRAEL}
\date{3.10.98}
\email{amnon@wisdom.weizmann.ac.il}

\begin{document}

\begin{abstract}
Let $k$ be a field and $A$ a noetherian (noncommutative)
$k$-algebra. The rigid dualizing complex of $A$ was introduced by
Van den Bergh. When $A = \opn{U}(\mfrak{g})$, the enveloping
algebra of a finite dimensional Lie algebra $\mfrak{g}$, Van den
Bergh conjectured that the rigid dualizing complex is
$(\opn{U}(\mfrak{g}) \otimes \bwedge^{n} \mfrak{g})[n]$,
where $n = \opn{dim} \mfrak{g}$.
We prove this conjecture, and give a few applications in
representation theory and Hochschild cohomology.
\end{abstract}

\maketitle

% ** section 0 **
\setcounter{section}{-1}
\section{Introduction}

{\em Dualizing complexes} were introduced as part of
Grothendieck Duality Theory on schemes, in \cite{RD}, and the
noncommutative version was first studied in \cite{Ye}.
The basic change is that a dualizing complex
over a noncommutative ring is a complex of bimodules. For
technical reasons we work with noetherian algebras over a base
field $k$, and abbreviate $\otimes := \otimes_{k}$.
Given an algebra $A$, we write $A^{\circ}$
for the opposite algebra, and
$A^{\mrm{e}} := A \otimes A^{\circ}$.
We consider left modules by default. A dualizing complex $R$
is an object in the bounded derived category of bimodules
$\cat{D}^{\mrm{b}}(\cat{Mod} A^{\mrm{e}})$,
of finite injective dimension on both sides,
such that the functors
$\mrm{R} \opn{Hom}_{A}(-, R)$
and
$\mrm{R} \opn{Hom}_{A^{\circ}}(-, R)$
induce a duality (i.e.\ a contravariant equivalence) between
$\cat{D}^{\mrm{b}}_{\mrm{f}}(\cat{Mod} A)$
and
$\cat{D}^{\mrm{b}}_{\mrm{f}}(\cat{Mod} A^{\circ})$.
The subscript f denotes complexes with finitely generated
cohomologies. See \cite{Ye} and \cite{YZ} for details on
noncommutative Grothendieck duality.

In the fundamental paper \cite{VdB1}, Van den Bergh defined the
{\em rigid dualizing complex} of a $k$-algebra $A$. A dualizing
complex $R$ is rigid if there exists an isomorphism
\begin{equation} \label{eqn0.1}
\rho: R \iso \mrm{R} \opn{Hom}_{A^{\mrm{e}}}(A, R \otimes R)
\end{equation}
in $\cat{D}(\cat{Mod} A^{\mrm{e}})$,
which we shall call a {\em rigidifying isomorphism}.
According to \cite{VdB1}, a rigid dualizing complex $R$, if it
exists, is unique up to isomorphism. Moreover it turns out that
rigid dualizing complexes are functorial with respect to
finite homomorphisms of $k$-algebras (under some technical
restrictions; cf.\ Theorem \ref{thm1.1}).

For instance, if $A$ is a commutative finite type $k$-algebra,
$\pi : X = \opn{Spec} A \ar \opn{Spec} k$ is the structural
morphism and
$\pi^{!} : \cat{D}^{\mrm{b}}_{\mrm{f}}(\cat{Mod} k) \ar
\cat{D}^{\mrm{b}}_{\mrm{f}}(\cat{Mod} A)$
is the twisted inverse image of \cite{RD}, then
$R := \pi^{!} k$ is a rigid dualizing complex, and $\rho$
is the fundamental class of the diagonal
$X \inj X \times X$.

Regarding existence of rigid dualizing complexes, Van den Bergh
proved the following result: if $A$ is
filtered such that $B := \opn{gr} A$ is a connected graded
noetherian $k$-algebra, and $B$ has a {\em balanced dualizing
complex} in the sense of \cite{Ye}, then $A$ has a rigid
dualizing complex.
In particular this holds for $A = \opn{U}(\mfrak{g})$,
the universal enveloping algebra of a finite dimensional Lie
algebra $\mfrak{g}$.

Our main result verifies a conjecture of Van den Bergh
(private communication, 1996):

\begin{thm} \label{thm0.1}
Let $\mfrak{g}$ be a finite dimensional Lie algebra over $k$. Then
the rigid dualizing complex of the universal enveloping algebra
$\opn{U}(\mfrak{g})$ is
\[ R =
\bigl( \mrm{U}(\mfrak{g}) \otimes \bwedge^{n} \mfrak{g}
\bigr)[n] , \]
where $n = \opn{dim} \mfrak{g}$, and we consider
$\bwedge^{n} \mfrak{g}$ as a $\mrm{U}(\mfrak{g})$-bimodule with
trivial action from the left and adjoint action from the right.
\end{thm}

Observe that in the two extreme cases -- $\mfrak{g}$ abelian or
semisimple -- the adjoint representation on
$\bwedge^{n} \mfrak{g}$ is trivial.
But for a solvable Lie algebra we can get something nontrivial,
as shown in Example \ref{exa1.1}.
The semisimple case was already known to Van den Bergh (cf.\
\cite{VdB2} Corollary 6).

An indication that Theorem \ref{thm0.1} should be true can be seen
by deforming $\mfrak{g}$ to an abelian Lie algebra. In the abelian
case $A = \mrm{U}(\mfrak{g})$ is a commutative polynomial
algebra, and there is a canonical isomorphism
$\mrm{U}(\mfrak{g}) \otimes \bwedge^{n} \mfrak{g}
\cong \Omega^{n}_{A / k}$.
As mentioned before, the complex
$\Omega^{n}_{A / k}[n] = \pi^{!} k$
is the rigid dualizing complex of $A$ (cf.\ Remark \ref{rem2.1}).

The proof of Theorem \ref{thm0.1} is at the end of
Section \ref{sec1}. In Section \ref{sec2}
we give a few corollaries of Theorem \ref{thm0.1}, and also an
analogous result for a ring $\mcal{D}(C)$ of differential operators
over a smooth commutative $k$-algebra $C$. 

\medskip \noindent \textbf{Acknowledgments.}\
I am grateful to Michel Van den Bergh for telling me about his
conjecture and for many helpful suggestions.
This paper was written during visits to MIT and the University of
Washington in 1998, and I wish to thank the Departments of
Mathematics at these universities for their hospitality, and
especially Michael Artin.

% ** section 1 **
\section{Proof of Main Result}
\label{sec1}

Let us start with a some general facts about rigid dualizing 
complexes of filtered $k$-algebras. 

If $\gamma$ is an automorphism of a ring $A$ then the
twist of a right module $M$ by $\gamma$ is $M_{\gamma}$, where the
new action is via $\gamma$.
In particular the twisted bimodule $A_{\gamma}$ has basis
$1_{\gamma}$, and
$1_{\gamma} \cdot a = \gamma(a) \cdot 1_{\gamma}$ for $a \in A$.
The shift by $i \in \mbb{Z}$ of a graded module $M$ is denoted by
$M(i)$, whereas the shift of a complex $M^{\bdot}$ is
$M^{\bdot}[i]$.

\begin{prop} \label{prop1.5}
Let $A$ be a filtered $k$-algebra, and assume $\opn{gr} A$ is a
connected graded, noetherian, Artin-Schelter Gorenstein algebra.
\begin{enumerate}
\item $A$ has a rigid dualizing complex
$R_{A} = \omega_{A}[n]$ for some integer $n$ and invertible
bimodule $\omega_{A}$.
Furthermore $\omega_{A} \cong A_{\gamma}$
where $\gamma$ is a filtered $k$-algebra automorphism of $A$.
\item The balanced dualizing complex of $\opn{gr} A$ is
$R_{\opn{gr} A} = \omega_{\opn{gr} A}[n]$, and
$\omega_{\opn{gr} A} \cong (\opn{gr} A)_{\opn{gr}(\gamma)}(m)$
for some integer $m$.
\end{enumerate}
\end{prop}

\begin{proof}
(Cf.\ \cite{YZ} Proposition 6.18.)
Let $\tilde{A} := \opn{Rees} A \subset A[t, t^{-1}]$
denote the Rees algebra. Recall that $t$ is a central variable
and  $(\opn{Rees} A)_{i} = F_{i} A \cdot t^{i}$.
Since $\tilde{A}$ is also AS-Gorenstein its balanced dualizing
complex is
$R_{\tilde{A}} = \tilde{A}_{\tilde{\gamma}}(m - 1)[n + 1]$
where $\tilde{\gamma}$ is a graded $k$-algebra automorphism
and $m, n \in \mbb{Z}$. Because
$\tilde{A}_{\tilde{\gamma}}$ is $k[t]$-central, $\tilde{\gamma}$ 
is in fact a $k[t]$-algebra automorphism.
Now by \cite{YZ} Theorem 6.2,
$R_{A} \cong
(\tilde{A}_{\tilde{\gamma}} \otimes_{\tilde{A}} A)[n]$.
On the other hand, using the exact sequence
$0 \ar \tilde{A}(-1) \xrightarrow{t} \tilde{A} \ar
\opn{gr} A \ar 0$
we get
\[ R_{\opn{gr} A} \cong
\opn{R} \opn{Hom}_{\tilde{A}}(\opn{gr} A,
\tilde{A}_{\tilde{\gamma}}(m - 1)[n + 1])
\cong (\tilde{A}_{\tilde{\gamma}} \otimes_{\tilde{A}} \opn{gr} A)
(m)[n] . \]
\end{proof}

We call $\omega_{A}$ the {\em dualizing bimodule} of $A$ and
$\gamma$ is the {\em dualizing automorphism}.

Next let us quote a result from \cite{YZ}.
A filtration $\{ F_{i} A \}$ is said to be noetherian connected
if $\opn{gr}^{F} A$ is a noetherian connected graded $k$-algebra. A 
ring homomorphism $A \ar B$ is finite centralizing if
$B = \sum_{i = 1}^{l} A \cdot b_{i}$
for some elements $b_{1}, \ldots, b_{l} \in B$
that commute with $A$. 

\begin{thm}[\tup{\cite{YZ}} Theorem 6.17] \label{thm1.1}
Let $A \ar B$ be a finite centralizing homomorphism of $k$-algebras.
Suppose $A$ has a noetherian connected filtration $\{ F_{i} A \}$
and $\opn{gr}^{F} A$ has a balanced dualizing complex. Then
the algebras $A$ and $B$ have rigid dualizing complexes
$R_{A}$ and $R_{B}$ respectively, and the trace morphism
$\opn{Tr}_{B / A} : R_{B} \ar R_{A}$ in
$\cat{D}(\cat{Mod} A^{\mrm{e}})$ exists. The trace induces
isomorphisms
\[  R_{B} \cong \opn{R} \opn{Hom}_{A}(B, R_{A}) \cong
\opn{R} \opn{Hom}_{A^{\circ}}(B, R_{A}) \]
in $\cat{D}(\cat{Mod} A^{\mrm{e}})$.
\end{thm}

Let $\mfrak{g}$ be a finite dimensional Lie algebra over the field
$k$, let $\mfrak{h} \subset \mfrak{g}$ be a subalgebra, and
denote by
$\mbf{K}_{\bdot}(\mfrak{h})$ the Chevalley-Eilenberg complex of
$\opn{U}(\mfrak{h})$, namely the free resolution of the trivial
$\mfrak{h}$-module $k$ (cf.\ \cite{CE} Section XIII.7 or \cite{Lo}
Section 10.1.3). Recall that for any $i$ one has
$\mbf{K}_{i}(\mfrak{h}) :=
\opn{U}(\mfrak{h}) \otimes \bwedge^{i} \mfrak{h}$,
a free left $\opn{U}(\mfrak{h})$-module (the action on the exterior
power $\bwedge^{i} \mfrak{h}$ is trivial). The boundary operator
$\delta : \mbf{K}_{i}(\mfrak{h}) \ar \mbf{K}_{i - 1}(\mfrak{h})$
is
\begin{multline*}
\delta(1 \otimes x_{1} \wedge \cdots \wedge x_{i}) =
\sum_{p = 1}^{i} (-1)^{p + 1} x_{p} \otimes
x_{1} \wedge \cdots \widehat{x}_{p} \cdots
\wedge x_{i} \\
+ \sum_{1 \leq p < q \leq i} (-1)^{p + q} \otimes [x_{p}, x_{q}]
\wedge x_{1} \wedge \cdots \widehat{x}_{p}
\cdots \widehat{x}_{q} \cdots \wedge x_{i}
\end{multline*}
for $x_{1}, \ldots, x_{i} \in \mfrak{h}$.
Define
\[ \mbf{K}_{i}(\mfrak{g}; \mfrak{h}) :=
\opn{U}(\mfrak{g}) \otimes_{\opn{U}(\mfrak{h})}
\mbf{K}_{i}(\mfrak{h}) \cong
\opn{U}(\mfrak{g}) \otimes \bwedge^{i} \mfrak{h} , \]
so that
$(\mbf{K}_{\bdot}(\mfrak{g}; \mfrak{h}), \delta)$
is a complex of free left $\opn{U}(\mfrak{g})$-modules.
As usual for any two $\opn{U}(\mfrak{g})$-modules $M, N$ the
tensor product $M \otimes N$ is also a $\opn{U}(\mfrak{g})$-module
by the coproduct.

\begin{lem} \label{lem1.2}
Suppose $\mfrak{h} \subset \mfrak{g}$ is an ideal, and consider
$\bwedge^{i} \mfrak{h}$ as a right $\opn{U}(\mfrak{g})$-module
by the adjoint action, so that
$\mbf{K}_{i}(\mfrak{g}; \mfrak{h})$
becomes a $\opn{U}(\mfrak{g})$-bimodule.
\begin{enumerate}
\item The boundary operator
$\delta: \mbf{K}_{i}(\mfrak{g}; \mfrak{h}) \ar
\mbf{K}_{i - 1}(\mfrak{g}; \mfrak{h})$
commutes with the right $\opn{U}(\mfrak{g})$-action.
\item There is a quasi-isomorphism of complexes of
$\opn{U}(\mfrak{g})$-bimodules
$\mbf{K}^{\bdot}(\mfrak{g}; \mfrak{h}) \ar
\opn{U}(\mfrak{g} / \mfrak{h})$.
\end{enumerate}
\end{lem}

\begin{proof}
1. Since
$\bwedge^{i} \mfrak{h} \subset \bwedge^{i} \mfrak{g}$
is a $\opn{U}(\mfrak{g})$-submodule for the adjoint action,
it follows that
$\mbf{K}_{i}(\mfrak{g}; \mfrak{h}) \subset
\mbf{K}_{i}(\mfrak{g})$
is a sub $\opn{U}(\mfrak{g})$-bimodule.
Hence we may assume that $\mfrak{h} = \mfrak{g}$
and
$\mbf{K}_{\bdot}(\mfrak{g}; \mfrak{h}) =
\mbf{K}_{\bdot}(\mfrak{g})$.
But then the assertion is \cite{Lo} Proposition 10.1.7.
(I wish to thank P. Smith for referring me to \cite{Lo}.)

\medskip \noindent
2. As usual we let
$\mbf{K}^{i}(\mfrak{g}; \mfrak{h}) :=
\mbf{K}_{-i}(\mfrak{g}; \mfrak{h})$, and the coboundary operator
is
$(-1)^{i + 1} \delta : \mbf{K}^{i}(\mfrak{g}; \mfrak{h}) \ar
\mbf{K}^{i + 1}(\mfrak{g}; \mfrak{h})$.
Since $\opn{U}(\mfrak{h}) \ar \opn{U}(\mfrak{g})$ is flat we get
$\mrm{H}^{i} \mbf{K}^{\bdot}(\mfrak{g}; \mfrak{h}) = 0$
if $i < 0$. For $i = 0$ we note that
$\opn{U}(\mfrak{g}) \cdot \mfrak{h} = \mfrak{h} \cdot
\opn{U}(\mfrak{g})$
is a two-sided ideal, and
\[ \opn{U}(\mfrak{g} / \mfrak{h}) \cong
\opn{U}(\mfrak{g})  / \opn{U}(\mfrak{g}) \cdot \mfrak{h} \cong
\mrm{H}^{0} \mbf{K}^{\bdot}(\mfrak{g}; \mfrak{h}) \]
as $\opn{U}(\mfrak{g})$-bimodules.
\end{proof}

For any $k$-module $M$ let 
$M^{*} := \opn{Hom}_{k}(M, k)$.
We consider $\bwedge^{n} \mfrak{g}^{*}$ as a right
$\opn{U}(\mfrak{g})$-module with the coadjoint action, and a left
$\opn{U}(\mfrak{g})$-module with the trivial action.

\begin{lem} \label{lem1.4}
Let $\mfrak{h} \subset \mfrak{g}$ be an ideal, with
$\opn{dim}_{k} \mfrak{h} = m$.
Assume that
$\gamma(\opn{U}(\mfrak{g}) \cdot \mfrak{h}) =
\opn{U}(\mfrak{g}) \cdot \mfrak{h}$.
Then
\[ \opn{Ext}^{q}_{\opn{U}(\mfrak{g})} \bigl(\opn{U}(\mfrak{g} /
\mfrak{h}), \opn{U}(\mfrak{g}) \bigr) \cong
\begin{cases}
\opn{U}(\mfrak{g} / \mfrak{h}) \otimes \bwedge^{m}
\mfrak{h}^{*} & \text{if } q = m \\
0 & \text{if } q \neq m
\end{cases}  \]
as $\opn{U}(\mfrak{g})$-bimodules.
\end{lem}

\begin{proof}
Since $\opn{gr} \opn{U}(\mfrak{g})$ is a commutative polynomial
algebra in $n$ variables we know that its balanced dualizing
complex is
$R_{\opn{gr} \opn{U}(\mfrak{g})} \cong
(\opn{gr} \opn{U}(\mfrak{g})(-n)[n]$.
Therefore by Proposition \ref{prop1.5} the rigid dualizing
complexes of
$\opn{U}(\mfrak{g})$ and $\opn{U}(\mfrak{g} / \mfrak{h})$
are
$R_{\opn{U}(\mfrak{g})} \cong \opn{U}(\mfrak{g})_{\gamma}[n]$ and
$R_{\opn{U}(\mfrak{g} / \mfrak{h})} \cong
\opn{U}(\mfrak{g} / \mfrak{h})_{\tau}[n - m]$,
respectively, where $\tau$ is the dualizing automorphism of
$\opn{U}(\mfrak{g} / \mfrak{h})$.
According to Theorem \ref{thm1.1} we get the vanishing of all
$\opn{Ext}^{q}$, $q \neq m$, and
\[ M := \opn{Ext}^{m}_{\opn{U}(\mfrak{g})}
\bigl(
\opn{U}(\mfrak{g} / \mfrak{h}), \opn{U}(\mfrak{g}) \bigr)
\cong \opn{U}(\mfrak{g} / \mfrak{h})_{\tau \gamma^{-1}} \]
as $\opn{U}(\mfrak{g})$-bimodules.

According to Lemma \ref{lem1.2} we get
\[ M =
\mrm{H}^{m} \opn{Hom}_{\opn{U}(\mfrak{g})}
\bigl(
\mbf{K}^{\bdot}(\mfrak{g}; \mfrak{h}),
\opn{U}(\mfrak{g}) \bigr) , \]
so the bimodule $M$ is a quotient of
$\opn{U}(\mfrak{g}) \otimes \bwedge^{m} \mfrak{h}^{*}$.
Let $\alpha$ be any $k$-basis of $\bwedge^{m} \mfrak{h}^{*}$,
and let $\beta$ be the image of
$1 \otimes \alpha \in
\opn{U}(\mfrak{g}) \otimes \bwedge^{m} \mfrak{h}^{*}$
in the $\opn{U}(\mfrak{g} / \mfrak{h})$-bimodule $M$.
Hence for any $x \in \mfrak{g}$ we have
\[ \beta \cdot x =
(x - \opn{tr}(\opn{ad}_{\bwedge^{m} \mfrak{h}^{*}} x))
\cdot \beta . \]
Since $M$ is free of rank $1$ on either side as
$\opn{U}(\mfrak{g} / \mfrak{h})$-module, and since
$\opn{U}(\mfrak{g} / \mfrak{h})$ is an integral domain, it follows
that the generator $\beta$ is a basis of $M$.
Sending
$\beta \mapsto 1 \otimes \alpha \in
\opn{U}(\mfrak{g} / \mfrak{h}) \otimes \bwedge^{m} \mfrak{h}^{*}$
is the desired isomorphism of $\opn{U}(\mfrak{g})$-bimodules.
\end{proof}

Here is another result of Van den Bergh (cf.\ \cite{VdB2}, proof
of Corollary 6).

\begin{lem} \label{lem1.6}
Let $A$ be a positively filtered $k$-algebra such that
$\opn{gr} A$ is commutative and $\opn{gr}_{0} A = k$.
Let $\mfrak{g} := \opn{gr}_{1} A$, so $\mfrak{g}$
is a Lie algebra over $k$.
Let $\gamma$ be a filtered $k$-algebra automorphism of $A$
such that $\opn{gr}(\gamma)$ is the identity. Then there is a
Lie homomorphism $\lambda : \mfrak{g} \ar k$
such that
$\gamma(a) = a + \lambda(\bar{a})$ for all $a \in F_{1} A$,
where $\bar{a} \in \mfrak{g}$ is the symbol of $a$.
\end{lem}

\begin{proof}
Define $\lambda(a) := \gamma(a) - a$ for $a \in F_{1} A$. It factors
through $F_{1} A \surj \mfrak{g} \ar F_{0} A \inj F_{1} A$,
is easily seen to be $k$-linear, and $\lambda([a, b]) = 0$.
\end{proof}

At last here is the proof of our main result.

\begin{proof}[Proof of Theorem \tup{\ref{thm0.1}}]
According to Proposition \ref{prop1.5}, the rigid dualizing
complex of $\opn{U}(\mfrak{g})$ is
$R_{\opn{U}(\mfrak{g})} \cong \opn{U}(\mfrak{g})_{\gamma}[n]$;
and $\opn{gr}(\gamma)$ is the identity. In view of Lemma
\ref{lem1.6}, it remains to prove that
$\lambda = - \opn{tr} \opn{ad}_{\bwedge^{n} \mfrak{g}}$.
Since $\lambda$ is a Lie homomorphism it has to vanish on
the commutator ideal
$\mfrak{h} := [\mfrak{g}, \mfrak{g}]$,
and so it factors through $\mfrak{a} := \mfrak{g} / \mfrak{h}$.
Therefore it suffices to prove that the induced automorphism
$\bar{\gamma}$ of $\opn{U}(\mfrak{a})$ satisfies
$\bar{\gamma}(y) = y
- \opn{tr} (\opn{ad}_{\bwedge^{n} \mfrak{g}} y)$
for $y \in \mfrak{a}$.

The algebra $\opn{U}(\mfrak{a})$ is a
commutative polynomial algebra in $l = n - m$ variables,
where $m = \opn{dim}_{k} \mfrak{h}$, so its rigid dualizing
complex is $\opn{U}(\mfrak{a})[l]$.
According to Lemma \ref{lem1.4} and Theorem \ref{thm1.1}
we get
\[ \opn{U}(\mfrak{a}) \cong
\opn{Ext}^{m}_{\opn{U}(\mfrak{g})} \bigl(\opn{U}(\mfrak{a}),
\opn{U}(\mfrak{g})_{\gamma} \bigr) \cong
\opn{U}(\mfrak{a})_{\gamma} \otimes \bwedge^{m} \mfrak{h}^{*} \]
as $\opn{U}(\mfrak{g})$-bimodules. Therefore
$\opn{U}(\mfrak{a})_{\bar{\gamma}} \cong
\opn{U}(\mfrak{a}) \otimes \bwedge^{m} \mfrak{h}$,
so
$\bar{\gamma}(y) = y -
\opn{tr}(\opn{ad}_{\bwedge^{m} \mfrak{h}} y)$
for all $y \in \mfrak{a}$. Finally, since
$\bwedge^{n - m} \mfrak{a}$ is a trivial representation of
 $\mfrak{g}$, one has
$\bwedge^{m} \mfrak{h} \cong \bwedge^{n} \mfrak{g}$.
\end{proof}

\begin{que} \label{que1.1}
Suppose $\mfrak{g}$ is semisimple and $\opn{char} k = 0$.
Does the quantum enveloping algebra $\opn{U}_{q}(\mfrak{g})$
admit a rigid dualizing complex? If so, what is it?
\end{que}

% ** section 2 **
\section{Some Corollaries and Complements}
\label{sec2}

\begin{cor} \label{cor2.3}
Let $M$ be any finitely generated $\opn{U}(\mfrak{g})$-module,
pure of $\opn{GKdim} = m$, and let
$I := \opn{Ann}_{\opn{U}(\mfrak{g})} M$. Then
\[ \opn{Ann}_{\opn{U}(\mfrak{g})^{\circ}}
\opn{Ext}^{n - m}_{\opn{U}(\mfrak{g})} 
\bigl( M, \opn{U}(\mfrak{g}) \bigr)
= \gamma(I) \subset \opn{U}(\mfrak{g})^{\circ} , \]
where $\gamma$ is the dualizing automorphism.
\end{cor}

\begin{proof}
Let us view $\gamma$ as an anti-isomorphism
$\gamma: \opn{U}(\mfrak{g}) \ar \opn{U}(\mfrak{g})^{\circ}$.
Define
$M':= \opn{Ext}^{n - m}_{\opn{U}(\mfrak{g})}
\bigl( M, \opn{U}(\mfrak{g}) \bigr)$
and
$I':= \opn{Ann}_{\opn{U}(\mfrak{g})^{\circ}} M'$.
By \cite{YZ} Proposition 6.18(4) one has
$\gamma(I) \subset I'$. Since $M$ is pure,
$M \subset M'' := \opn{Ext}^{n - m}_{\opn{U}( \mfrak{g})^{\circ}} 
\bigl( M', \opn{U}(\mfrak{g}) \bigr)$.
Hence
$\gamma^{-1}(I') \subset \opn{Ann}_{\opn{U}(\mfrak{g})} M''
\subset I$.
\end{proof}

It is a standard fact that if $M$ is a finite dimensional
representation of $\mfrak{g}$, then
$\opn{Ext}^{q}_{\opn{U}(\mfrak{g})} 
\bigl( M, \opn{U}(\mfrak{g}) \bigr) = 0$
for $q < n$. The group
$\opn{Ext}^{n}_{\opn{U}(\mfrak{g})} 
\bigl( M, \opn{U}(\mfrak{g}) \bigr)$ 
is a right $\opn{U}(\mfrak{g})$-module, but the structure is
not obvious. Since we can make $M$ into a
$\opn{U}(\mfrak{g})$-bimodule with trivial
right action, the next corollary gives the answer.

\begin{cor} \label{cor2.2}
Suppose $M$ is a finite dimensional $k$-central
$\opn{U}(\mfrak{g})$-bimodule. Then there is an isomorphism
of $\opn{U}(\mfrak{g})$-bimodules
\[ \opn{Ext}^{n}_{\opn{U}(\mfrak{g})}
\bigl( M, \opn{U}(\mfrak{g}) \bigr) \cong
M^{*} \otimes \bwedge^{n} \mfrak{g}^{*} , \]
which is functorial in $M$.
\end{cor}

\begin{proof}
Let
$I := \opn{Ann}_{\opn{U}(\mfrak{g})} M$
and $B := \opn{U}(\mfrak{g}) / I$. Since $k \ar B$ is a finite
homomorphism the rigid dualizing complex of $B$ is
$B^{*} = \opn{Hom}_{k}(B, k)$. By \cite{YZ} Proposition 3.9,
\[  \opn{Ext}^{n}_{\opn{U}(\mfrak{g})}
\bigl( M, \opn{U}(\mfrak{g}) \otimes \bwedge^{n} \mfrak{g}
\bigr) \cong
\opn{Hom}_{B}(M, B^{*}) \cong M^{*}  \]
as $\opn{U}(\mfrak{g})$-bimodules. Now twist by
$\bwedge^{n} \mfrak{g}^{*}$.
\end{proof}

Theorem \ref{thm0.1} has an interpretation in terms of Hochschild
cohomology. For a $\opn{U}(\mfrak{g})$-bimodule $M$ denote by
$\opn{H}^{q}(\opn{U}(\mfrak{g}), M)$ and
$\opn{H}_{q}(\opn{U}(\mfrak{g}), M)$ the
Hochschild cohomology and homology, respectively.

\begin{cor} \label{cor2.5}
There are $\opn{U}(\mfrak{g})$-bimodule isomorphisms
\[ \opn{H}^{q} \bigl( \opn{U}(\mfrak{g}), 
\opn{U}(\mfrak{g})^{\mrm{e}} \bigr) \cong
\begin{cases}
\opn{U}(\mfrak{g}) \otimes \bwedge^{n} \mfrak{g}^{*}
& \tup{ if } q = n \\
0 & \tup{ if } q \neq n .
\end{cases} \]
\end{cor}

\begin{proof}
Let's write
$\omega := \omega_{\opn{U}(\mfrak{g})}$
and
$\omega^{\vee} := \opn{Hom}_{\opn{U}(\mfrak{g})}(\omega,
\opn{U}(\mfrak{g}))$.
By formula (\ref{eqn0.1}),
$\omega \cong
\opn{Ext}^{n}_{\opn{U}(\mfrak{g})} \bigl( \opn{U}(\mfrak{g}),
\omega \otimes \omega \bigr)$
as bimodules, so applying the twist
$- \otimes_{\opn{U}(\mfrak{g})^{\mrm{e}}}
(\omega^{\vee} \otimes \omega^{\vee})$
we get
$\omega^{\vee} \cong
\opn{Ext}^{n}_{\opn{U}(\mfrak{g})} \bigl( \opn{U}(\mfrak{g}),
\opn{U}(\mfrak{g})^{\mrm{e}} \bigr)$.
But by Theorem \ref{thm0.1},
$\omega^{\vee} \cong
\opn{U}(\mfrak{g}) \otimes \bwedge^{n} \mfrak{g}^{*}$.
\end{proof}

In \cite{VdB2}, Van den Bergh proves a Poincar\'{e} duality
between the Hochschild cohomology and homology of certain
Gorenstein algebras $A$. We obtain the following variation
of his result.

\begin{cor} \label{cor2.4}
Let $M$ be any $k$-central $\opn{U}(\mfrak{g})$-bimodule. Then
\[ \opn{H}^{q} \bigl( \opn{U}(\mfrak{g}), M \bigr) \cong
\opn{H}_{n - q} \bigl( \opn{U}(\mfrak{g}),
M \otimes \bwedge^{n} \mfrak{g}^{*} \bigr) . \]
\end{cor}

\begin{proof}
Corollary \ref{cor2.5} says that 
\[ \opn{R} \opn{Hom}_{\opn{U}(\mfrak{g})^{\mrm{e}}}
\bigl( \opn{U}(\mfrak{g}), \opn{U}(\mfrak{g})^{\mrm{e}} \bigr)[n]
\cong \omega^{\vee} \cong
\opn{U}(\mfrak{g}) \otimes \bwedge^{n} \mfrak{g}^{*} \]
in
$\cat{D}(\cat{Mod} \opn{U}(\mfrak{g})^{\mrm{e}})$.
Copying the proof of \cite{VdB2} Theorem 1 we obtain
\[ \begin{aligned}
\opn{H}^{q} \bigl( \opn{U}(\mfrak{g}), M \bigr) & \cong
\opn{H}^{q} \opn{R} \opn{Hom}_{\opn{U}(\mfrak{g})^{\mrm{e}}}
\bigl( \opn{U}(\mfrak{g}), M) \\
& \cong \opn{H}^{q} \Bigl( 
\opn{R} \opn{Hom}_{\opn{U}(\mfrak{g})^{\mrm{e}}}
\bigl( \opn{U}(\mfrak{g}), \opn{U}(\mfrak{g})^{\mrm{e}} \bigr)
\otimes^{\mrm{L}}_{\opn{U}(\mfrak{g})^{\mrm{e}}} M \Bigr) \\
& \cong \opn{H}^{q - n} \bigl(\omega^{\vee} 
\otimes^{\mrm{L}}_{\opn{U}(\mfrak{g})^{\mrm{e}}} M) \\
& \cong \opn{H}^{q - n} \bigl( \opn{U}(\mfrak{g})
\otimes^{\mrm{L}}_{\opn{U}(\mfrak{g})^{\mrm{e}}}
( M \otimes_{\opn{U}(\mfrak{g})} \omega^{\vee} )  \bigr) \\
& \cong \opn{H}_{n - q} \bigl( \opn{U}(\mfrak{g}),
M \otimes \bwedge^{n} \mfrak{g}^{*} \bigr) .
\end{aligned} \]
\end{proof}

Here is an easy example where the dualizing bimodule 
$\omega$ is not trivial.

\begin{exa} \label{exa1.1}
Let $\mfrak{g}$ be the nonabelian $2$-dimensional Lie algebra,
with basis $x, y$ such that $[x, y] = y$. Then
$\opn{tr}(\opn{ad}_{\bwedge^{2} \mfrak{g}} x) = 1$.
\end{exa}

If $\opn{char} k = 0$ and $C$ is a smooth, integral, commutative 
$k$-algebra then the ring of differential operators $\mcal{D}(C)$ is 
noetherian and has finite global dimension. Since $\mcal{D}(C)$ 
can be deformed to a smooth commutative $k$-algebra (namely the 
algebra of functions on the cotangent bundle of $\opn{Spec} C$), 
one could expect $\mcal{D}(C)$ to have a rigid dualizing 
complex. This is indeed true, and follows from results in 
$\mcal{D}$-module theory.

\begin{thm} \label{thm2.1}
Let $C$ be a smooth, integral, commutative $k$-algebra of dimension 
$n$, and assume  $\opn{char} k = 0$. 
Let $\mcal{D}(C)$ be the ring of differential operators.
Then the rigid dualizing complex of $\mcal{D}(C)$ is
$\mcal{D}(C)[2n]$.
\end{thm}

\begin{proof}
Let $X := \opn{Spec} C$ and
$X^{\mrm{e}} := X \times X \cong \opn{Spec} C^{\mrm{e}}$.
Then 
$\Gamma(X, \mcal{D}_{X}) \cong \mcal{D}(C)$,
$\Gamma(X^{\mrm{e}}, \mcal{D}_{X^{\mrm{e}}}) \cong 
\mcal{D}(C) \otimes \mcal{D}(C)$
and
$\mcal{D}(C)^{\circ} \cong \omega_{C} \otimes_{C} \mcal{D}(C)
\otimes_{C} \omega_{C}^{\vee}$. 

The sheaf 
$\mcal{D}_{X} \otimes_{\mcal{O}_{X}} \omega_{X}^{\vee}$
is filtered, and has two commuting left $\mcal{D}_{X}$-module 
structures. The two structures coincide on
$\opn{gr} (\mcal{D}_{X} \otimes_{\mcal{O}_{X}} \omega_{X}^{\vee})
\cong (\opn{gr} \mcal{D}_{X}) \otimes_{\mcal{O}_{X}}
\omega_{X}^{\vee}$. Hence there is an involution of 
$\mcal{D}_{X} \otimes_{\mcal{O}_{X}} \omega_{X}^{\vee}$,
which is the identity on the subsheaf 
$\omega_{X}^{\vee} = F_{0}(\mcal{D}_{X} \otimes_{\mcal{O}_{X}}
\omega_{X}^{\vee})$,
and exchanges the two $\mcal{D}_{X}$-module structures.

Denote by $\mbf{D}_{X}$ the duality functor on left 
$\mcal{D}_{X}$-modules, namely 
$\mbf{D}_{X} \mcal{M} :=$ \linebreak
$\opn{R} \mcal{H}om_{\mcal{D}_{X}}(\mcal{M}, 
\mcal{D}_{X} \otimes_{\mcal{O}_{X}} \omega_{X}^{\vee})[n]$;
cf.\ \cite{Bo} VI.3.6. 
Let $f : X \inj X^{\mrm{e}}$ be the diagonal embedding.
According to \cite{Bo} Proposition VII.9.6 there is a functorial 
isomorphism
$\mbf{D}_{X^{\mrm{e}}} \: f_{+} \cong f_{+} \: \mbf{D}_{X}$.
We shall apply this isomorphism with the $\mcal{D}_{X}$-module
$\mcal{O}_{X}$.

First note that $\mbf{D}_{X}\: \mcal{O}_{X} \cong \mcal{O}_{X}$ ,
as can be checked using the quasi-isomorphism 
$\Omega_{X}^{\bdot}(\mcal{D}_{X})[n] \otimes_{\mcal{O}_{X}}
\omega_{X}^{\vee} \ar \mcal{O}_{X}$
in $\cat{Mod} \mcal{D}_{X}$; cf.\ \cite{Bo} VI.3.5.
Next, by \cite{Bo} Theorem VI.7.4(ii) and Theorem VI.7.11 
(Kashiwara's Theorem) we see that
$f_{+}\: \mcal{O}_{X} \cong \mcal{D}_{X} \otimes_{\mcal{O}_{X}}
\omega_{X}^{\vee}$
in $\cat{Mod} \mcal{D}_{X^{\mrm{e}}}$.
Thus we have an isomorphism of $\mcal{D}_{X^{\mrm{e}}}$-modules
\[ \mcal{D}_{X} \otimes_{\mcal{O}_{X}} \omega_{X}^{\vee}
\cong
\mcal{E}xt^{2n}_{\mcal{D}_{X^{\mrm{e}}}} \bigl(
\mcal{D}_{X} \otimes_{\mcal{O}_{X}} \omega_{X}^{\vee},
\mcal{D}_{X^{\mrm{e}}} \otimes_{\mcal{O}_{X^{\mrm{e}}}} 
\omega_{X^{\mrm{e}}}^{\vee} \bigr) . \]
Passing to global sections, replacing $\mcal{D}(C)$ by
$\mcal{D}(C)^{\circ}$ and using the involution of
$\mcal{D}(C) \otimes_{C} \omega_{C}^{\vee}$, 
we get
\[ \begin{split}
& \mcal{D}(C) \otimes_{C} \omega_{C}^{\vee} \\[0.5ex]
& \qquad \cong
\opn{Ext}^{2n}_{\mcal{D}(C) \otimes \mcal{D}(C)}
\bigl( \mcal{D}(C) \otimes_{C} \omega_{C}^{\vee}, 
(\mcal{D}(C) \otimes_{C} \omega_{C}^{\vee}) \otimes
(\mcal{D}(C) \otimes_{C} \omega_{C}^{\vee}) \bigr) \\[0.5ex]
& \qquad \cong
\opn{Ext}^{2n}_{\mcal{D}(C) \otimes \mcal{D}(C)^{\circ}}
\bigl( \mcal{D}(C), 
(\mcal{D}(C) \otimes_{C} \omega_{C}^{\vee}) \otimes
\mcal{D}(C) \bigr) \\[0.5ex]
& \qquad \cong
\opn{Ext}^{2n}_{\mcal{D}(C)^{\mrm{e}}}
\bigl( \mcal{D}(C), 
\mcal{D}(C) \otimes \mcal{D}(C) \bigr) 
\otimes_{C} \omega_{C}^{\vee} .
\end{split} \]
Twisting by $\omega_{C}$ and shifting degrees we obtain
an isomorphism
\[ \mcal{D}(C)[2n] \cong
\opn{R} \opn{Hom}_{\mcal{D}(C)^{\mrm{e}}}
\bigl( \mcal{D}(C), 
\mcal{D}(C)[2n] \otimes \mcal{D}(C)[2n] \bigr) \]
in 
$\cat{D}(\cat{Mod} \mcal{D}(C)^{\mrm{e}})$.
\end{proof}

By the same arguments given for Corollaries \ref{cor2.5} and
\ref{cor2.4}, one has:

\begin{cor} \label{cor2.6}
Let $\mcal{D}(C)$ be as above. Then there are
$\mcal{D}(C)$-bimodule isomorphisms
\[ \opn{H}^{q} \bigl( \mcal{D}(C), \mcal{D}(C)^{\mrm{e}} \bigr)
\cong
\begin{cases}
\mcal{D}(C)
& \tup{ if } q = 2n \\
0 & \tup{ if } q \neq 2n .
\end{cases} \]
For any $k$-central $\mcal{D}(C)$-bimodule $M$ one has
\[ \opn{H}^{q} \bigl( \mcal{D}(C), M \bigr) \cong
\opn{H}_{2n - q} \bigl( \mcal{D}(C), M \bigr) . \]
\end{cor}

\begin{rem} \label{rem2.1}
One can show that there is a canonical choice for the
rigidifying isomorphism $\rho$ of the complex
$R = \omega[n]$,
$\omega = \opn{U}(\mfrak{g}) \otimes \bwedge^{n} \mfrak{g}$.
This amounts to choosing an isomorphism of bimodules
$\rho : \omega \cong E^{n}(\opn{U}(\mfrak{g}))$, where
$E^{n}(\opn{U}(\mfrak{g})) :=
\opn{Ext}^{n}_{\opn{U}(\mfrak{g})^{\mrm{e}}}
\bigl( \opn{U}(\mfrak{g}), \omega \otimes \omega \bigr)$.
Here is a sketch of the proof. Let
$A := \opn{gr} \opn{U}(\mfrak{g}) = \opn{S}(\mfrak{g})$.
The bimodule $\omega$ is filtered, and there is a canonical 
isomorphism
$\opn{gr} \omega \cong \Omega^{n}_{A / k}$.
The standard spectral sequence of the filtration identifies
$\opn{gr} E^{n}(\opn{U}(\mfrak{g}))$ with
$E^{n}(A) :=
\opn{Ext}^{n}_{A^{\mrm{e}}}(A, \Omega^{2n}_{A^{\mrm{e}} / k})$.
But as mentioned in the Introduction,
$\Omega^{n}_{A / k}$ is the rigid dualizing complex of $A$, and it
comes equipped with a canonical isomorphism
$\Omega^{n}_{A / k} \iso E^{n}(A)$.
This isomorphism determines $\rho$.
A similar statement holds for Theorem \ref{thm2.1}.
As a consequence the isomorphisms of Corollaries \ref{cor2.5},
\ref{cor2.4} and \ref{cor2.6} are canonical.
(I thank Van den Bergh for mentioning this idea to me.)
\end{rem}

\end{document}